\documentstyle[11pt]{article}
\newcommand\be{\begin{equation}}
\newcommand\ee{\end{equation}}
\newtheorem{theorem}{Theorem}
\newtheorem{lemma}{Lemma}
\newtheorem{defn}{Definition}
\newtheorem{cor}{Corollary}
\def\ra{\rightarrow}
\def\P{{\cal P}}

\def\L{{\cal L}}
\def\A{{\cal A}}
\def\bA{{\bf A}}

\def\t*{\theta_{G^*}}
\def\m{{\bf m}}
\def\a{\alpha}
\def\n{{\bf n}}
\def\o{\omega}
\def\to{\tilde{\omega}}
\def\N{{\bf N}}
\def\L{{\cal L}}
\def\g{{\scriptstyle {\cal G}}}
\def\h{{\scriptstyle {\cal H}}}
\def\Otu{\Omega_{(t, u)}}
\begin{document}
\begin{titlepage}
\begin{center}
{\Large\bf  On Poisson actions of compact Lie groups}\\
\vskip 0.3cm

{\Large\bf  on symplectic manifolds} \\
\vskip 1cm
{\bf A. Yu. Alekseev\footnote{On leave of absence from Steklov
Institute, St.-Petersburg}}\\
\vskip 0.2cm
{\it Institute for Theoretical Physics, Uppsala University, Box
803}\\
{\it S-75108, Uppsala, Sweden\footnote{Permanent address}}\\
\vskip 0.1cm
{\it and}
\vskip 0.1cm
{\it Institut f\"ur Theoretische Physik, ETH, CH-8093 Z\"urich,
Switzerland}\\
\vskip 1.2cm
February 1996
\vskip 1.2cm
\end{center}

\begin{abstract}

Let  $G_{\P}$ be a compact simple Poisson-Lie group equipped with a
Poisson structure $\P$ and $(M, \o)$ be a symplectic manifold. Assume
that $M$ carries a Poisson action of $G_{\P}$ and there is an
equivariant
moment map in the sense of Lu and Weinstein which acts
to the dual Poisson-Lie group $G^*_{\P}$, $\m: M\rightarrow
G^*_{\P}$.  We prove
that  $M$ always possesses another symplectic form $\to$ so that the
$G$-action preserves $\tilde{\o}$ and there is a new moment map
$\mu= e^{-1} \circ \m: M\rightarrow \g^*$.
Here  $e$ is a universal (independent of $M$) invertible equivariant
map
$e: \g^*\rightarrow G^*_{\P}$. We suggest    new short proves of the
convexity theorem for the Poisson-Lie moment map, Poisson reduction
theorem
and the Ginzburg-Weinstein theorem on the isomorphism of $\g^*$ and
$G^*_{\P}$ as
Poisson spaces.

\end{abstract}
\end{titlepage}

The main goal of this paper is to compare Hamiltonian and Poisson
actions of
compact simple Lie groups on symplectic manifolds.  We prove that one
can always
exchange the Poisson action to the Hamiltonian one by an appropriate
change of
the symplectic structure. This trick reduces many questions
concerning Poisson actions
to their well known counterparts from the theory of Hamiltonian
$G$-actions.  In particular,
we suggest   new simple proves of the convexity theorem for the
Poisson-Lie moment map \cite{FR}, Poisson reduction theorem \cite{Lu}
and the Ginzburg-Weinstein theorem  \cite{GW}.
The results of this paper were announced in \cite{AA}.

\subsection*{Compact Poisson-Lie groups}

\begin{defn}
Let $G$ be a simple connected simply connected compact  Lie group
and $\P$  be a Poisson bracket on $G$. This pair defines a
Poisson-Lie
group if   the multiplication map $G\times G\rightarrow G$ is a
Poisson
map.
\end{defn}

Up to a scalar factor  Poisson-Lie structures on $G$ are in one to
one
correspondence with Manin triples $(d, \g, \g^*)$.

\begin{defn}
A triple of Lie algebras $(d, \g, \g^*)$ is called a Manin
triple if  $d$ has an invariant nondegenerate bilinear form $k$ and
$\g$ and
$\g^*$ are maximal isotropic subalgebras of $d$ which together span
$d$:
\be
k(\g, \g)=k(\g^*, \g^*)=0.
\ee
\end{defn}

The algebra $d$ is also called a Drinfeld double of $\g$.
In our particular example $d$ always coincides with the complex Lie
algebra $\g^{C}$
considered as an algebra over  real numbers. The scalar product
$k$ is given by the imaginary part of the Killing form $K$ on
$\g^{C}$:
\be \label{sp}
k(a, b)= Im \ K (a, b).
\ee

Up to an isomorphism the Manin triples including $d$ and $\g$ are
classified by
real valued antisymmetric bilinear forms  on the Cartan subalgebra
$\h$ of $\g$ \cite{LS}:
\be \label{mt}
d=\g + \g^*_u.
\ee
For each such a form $u$, the dual Lie algebra $\g^*_u$ is defined as
a semi direct sum of
two subalgebras
\be
\g^*_u= \n+ \h^*_u.
\ee
Here $\n$ is the maximal nilpotent subalgebra in $\g^{C}$. We can
always assume that it is generated by all positive roots of $\g^{C}$.
 The other subspace $\h^*_u$ is defined as follows
\be
\h^*_u=\{ i(a+iu(a)), a\in \h \}
\ee
where  $u(a)$ is  a map from $\h$ to itself corresponding to  the
form $u$.
Antisymmetry of $u$ implies
\be
K(a, u(b))+K(u(a), b)=0.
\ee

Let us denote a Poisson structure corresponding to the Manin triple
(\ref{mt}) by $\P_u$. Rescaling this Poisson bracket by a real factor
$t$ we get a family  parametrized by pairs $(t, u)$:
\be \label{fam}
\P_{(t, u)}= t \P_u.
\ee
This family  provides a complete classification of Poisson structures
on  compact simple Lie groups (up to isomorphisms) if we add the set
of points
\be
(t\rightarrow 0, u=\frac{w}{t}) \ , \ w=const
\ee
parametrized by $w$ and lying at infinity of the space of parameters
of the family (\ref{fam}).
We shall refer to this special family as to the case of $t=0$. In the
main part of the paper we always assume that $t\neq 0$ and collect
some details on the case of $t=0$ in Appendix.

Let us remark that the Lie algebras $\g$ and $\g^*$ enter the picture
in a symmetric way.
This means that the connected simply connected group $G^*_u$
corresponding to
the Lie algebra $\g^*_u$ also carries a Poisson-Lie structure defined
by the Manin triple.

In our example the group $G^*_u$ is a semi-direct product of the
maximal nilpotent group
$\N$ in $G^C$ and the subgroup $H^*_u$ of the complexefication of the
Cartan torus
\be
H^*_u=\{ exp(a), a\in \h^*_u \}.
\ee
In particular, for $G=SU(N)$ and $u=0$ the group $H^*_0$ is formed by
diagonal
matrices of unit determinant with positive eigenvalues.
The elements of $G^*_u$ may be visualized by embedding into $G^{C}$:
\be
G^*=\{ N \ exp\{i(a+iu(a))\}, N\in \N, a\in \h\}.
\ee

Let $a\rightarrow \bar{a}$ be an anti-involution of $G^{C}$ which
singles
out the compact form. It is convenient to introduce a map
\be
f: a\rightarrow  a\bar{a}
\ee
which maps $G^*_u$ into a certain subspace $SG$ of  $G^{C}$
\be
SG=\{ exp\{ ia\}, a\in \g\}.
\ee
Observe that though the dual group $G^*_u$ depends on the choice of
$u$, the target of the map $f$ is always the same space $SG$.

There is another way to characterize $SG$:
\be
SG=\{ x\in G^{C}, \bar{x}=x\}.
\ee
The bar operation being {\em anti}-involution, $SG$ is not a group.
Using the fact that any element of $SG$ may be brought to the maximal
torus
by conjugation with some element of $G$, the Iwasawa decomposition
and
the uniqueness of a positive square root of a positive real number
one easily proves
that the map $f$ is in fact invertible. Let us define the following
map $e_{(t, u)}$
from $\g^*$ to $G^*_u$:
\be \label{et}
e_{(t, u)}= f^{-1} \circ j  \ , \  j = E \circ K=exp\{2it\ \cdot \}
\circ K.
\ee
Here  $K$ stays for the Killing form which converts $\g^*$ to $\g$,
the exponential
map $E$ with additional $i$ maps $\g$ to $SG$ and the last map
$f^{-1}$ identifies
$SG$ with $G^*$. Let $a$ be an element of $\g^*$ and $A=e_{(t,
u)}(a)$. Then
the definition (\ref{et}) implies
\be
\bA\equiv A\bar{A}=j(a)=exp\{ 2it K(a)\} .
\ee

Both spaces $\g^*$ and $G^*_u$ carry natural actions of the group
$G$.
The dual space to the Lie algebra carries the coadjoint action
$Ad^*$:
\be
K(Ad^*(g) a)= gK(a) g^{-1}.
\ee

The $G$-action on the group $G^*_u$  is defined by using a somewhat
generalized version of the  Iwasawa decomposition:
\be
g\cdot A= A^g\cdot g'.
\ee
This is an equality  in $G^{C}$. In the right hand side $g'\in G$ and
$a^g\in G^*$. Existence and uniqueness of $A^g$ and $g'$ are ensured
by the corresponding properties of the Iwasawa decomposition. For
historical reasons this action of $G$ on $G^*$ is called {\em
dressing action} \cite{Sem}. To make notations closer to the case of
$\g^*$ we sometimes denote
\be
A^g= AD^*(g) A.
\ee
Observe that
\be
\bA^g=A^g \bar{A^g}= gA\bar{A}g^{-1}=g\bA g^{-1}.
\ee
This simple observation proves the following lemma.

\begin{lemma}
 The map $e_{(t, u)}$  intertwines coadjoint and dressing actions of
$G$ on $\g^*$ and $G^*_u$:
\be
AD^*(g) e_{(t, u)}(a)= e_{(t, u)} (Ad^*(g) a).
\ee
\end{lemma}

\vskip 0.5cm

The map $e_{(t, u)}$ has been introduced in \cite{FR}. We shall
discuss some new properties of this map in the next sections.

\subsection*{Moment map in the sense of Lu and Weinstein}

Let us recall the definitions of  the moment map for Hamiltonian and
Poisson
group actions on symplectic manifolds.

Let $M$ be a symplectic manifold equipped with an action $\A$
 of a compact Lie group $G$:
\be
\A: G\times M\rightarrow M, \  \A(g, x)= x^g.
\ee
One can introduce a universal vector field $v$ taking values in the
space $\g^*$  so that for any element $\alpha \in \g$ there is a
vector field
\be
v_{\epsilon}=<v, \epsilon>=\A_* (\epsilon).
\ee
In the right hand side we treat $\alpha$ as a right invariant vector
field on $G$.

\begin{defn}
The action $\A$ is called Hamiltonian if it preserves the Poisson
structure
on $M$:
\be \label{ham}
\A_*(\P_M)=\P_M.
\ee
\end{defn}
The Poisson tensor $\P$ is assumed to be the inverse of the matrix of
the symplectic form
$\o$ on $M$.

We are specifically interested in  symplectic manifolds equipped with
the $G$-action
and an equivariant moment map.

\begin{defn}
The map $\mu: M\rightarrow \g^*$ is called a moment
map if it satisfies the following property:
\be \label{mom}
\o(\cdot,  v)=  \mu^* (da).
\ee
Here $da$ is the natural linear $1$-form on $\g^*$ taking values in
$\g^*$.
\end{defn}

Existence of
the moment map ensures the invariance of the symplectic form with
respect to the $G$-action.

\begin{defn}
The moment map is said to be equivariant if
\be
Ad^*(g) \mu(x)= \mu (x^g).
\ee
\end{defn}

Let  $(G, \P_G)$ be a compact Poisson-Lie group, the Poisson
structure $\P_G$ being one of the standard list  parametrized by
pairs  $(t, u)$ (see the previous section).

\begin{defn}
The action of $\A: G\times M\rightarrow M$ is called
a Poisson action if it preserves the  Poisson structure in the
following sense:
\be
\A_*(\P_G+\P_M)=\P_M.
\ee
\end{defn}

Notice the difference with the standard definition  (\ref{ham}). If
$M$ is equipped with a Poisson action of $G$, the symplectic
structure on $M$  {\em is not} invariant with respect to the
$G$-action.

A Poisson counterpart of the  notion of the moment map has been
defined in \cite{LW}.  \begin{defn}
Let $G$ be a compact Poisson-Lie group equipped with a Poisson
structure $\P_{(t, u)}$. Let $\A: G\times M \rightarrow M$ be a
Poisson action of $G$ on the symplectic manifold $M$
The map $\m: M\rightarrow G^*_u$  is called a moment map in the sense
of Lu and Weinstein if
\be
\o(\cdot, v)= \frac{1}{t}\m^* (dAA^{-1}),
\ee
where  $dAA^{-1}$ is a right-invariant Maurer-Cartan form on $G^*_u$.
\end{defn}

The equivariance condition for the Poisson moment map $\m$ looks as
follows:
\be
AD^*(g) (\m (x))= \m (x^g).
\ee

\vskip 0.5cm

\subsection*{Comparing Hamiltonian and Poisson actions}

Here we formulate and  prove the main result of the paper.

\begin{theorem}
Let $(M, \o)$ be a symplectic manifold which carries an action $\A$
of a compact Poisson-Lie group $G$ equipped with a Poisson bracket
$\P_{(t, u)}$. Assume that there exists an equivariant moment map
$\m: M\ra G^*$.  Then one can define another symplectic form $\to$ on
$M$ with the following properties:

1) $\to$ is preserved by $\A$;

2) $\to$ belongs to the same cohomology class as $\o$;

3) the map  $\mu = e_{(t. u)}^{-1}\circ \m$ provides  an equivariant
moment map for the $G$-action $\A$ with respect to the symplectic
structure $\to$.
\end{theorem}

The main technical tool for proving this theorem is provided by the
following lemma.

\begin{lemma}
There exists such a 2-form $\Otu$ on $\g^*$, so that the following
two properties are fulfilled:

\vskip 0.2cm
1) The form $\Otu$ is closed $d \Otu =0$.

\vskip 0.2cm
2) $\Otu(\cdot, v)= \frac{1}{t} e_{(t, u)}^* dA A^{-1} - da$.

\vskip 0.2cm
Here $v$ is   the universal vector field corresponding to the
coadjoint action of $G$ on $\g^*$,
$a\in g^*$ and  $A=e_{(t, u)} (a) \in G^*_u$.
\end{lemma}

{\sc Proof of  Lemma}.
It is convenient to introduce a special notation for $\a= K(a)\in
\g$.
Let us consider the following  2-form on $\g^*$:
\be \label{Otu}
\Otu=\frac{1}{4it}
\{ K^* \sum_{k=2}^{\infty} \frac{(2it)^k}{k!} K(ad^{k-2}(\a)
d\a\wedge d\a) +  e_{(t, u)}^* K(A^{-1}dA\wedge d\bar{A}
\bar{A}^{-1})\} .
\ee
We claim that it satisfies both conditions of  Lemma 2.

It is convenient to split $\Otu$ into two pieces:
\be
\Otu= \o_1+\o_2,
\ee
where
\begin{eqnarray}
\o_1=\frac{1}{4it} K^* \sum_{k=2}^{\infty} \frac{(2it)^k}{k!}
K(ad^{k-2}(\a) d\a\wedge d\a), \nonumber \\
\o_2 =\frac{1}{4it} e_{(t, u)}^* K(A^{-1}dA\wedge d\bar{A}
\bar{A}^{-1}).
\end{eqnarray}

1) A direct calculation shows
\begin{eqnarray} \label{do1}
d \o_2=\frac{1}{4it} e_{(t, u)}^* d\{ K(A^{-1}dA\wedge d\bar{A}
\bar{A}^{-1})\} = \nonumber  \\
-\frac{1}{4it} e_{(t, u)}^* \{ K( (A^{-1}dA)^2\wedge d\bar{A}
\bar{A}^{-1}) +K( (A^{-1}dA)\wedge (d\bar{A} \bar{A}^{-1})^2)=\\
-\frac{1}{12it} j^* K(d\bA \bA^{-1}\wedge (d\bA \bA^{-1})^2).
\nonumber
\end{eqnarray}
Let us recall that $\bA=A\bar{A}=j(a)$.

Using equation
\be \label{dA}
d\bA \bA^{-1}=(E^{-1})^* \left( \frac{e^{2it
\lambda}-1}{\lambda}\right)_{\lambda=ad(a)} d\a
\ee
one can easily show  that
\begin{eqnarray} \label{do2}
d\o_1=d \{ \frac{1}{4it} K^* \sum_{k=2}^{\infty} \frac{(2it)^k}{k!}
K(ad^{k-2}(\a) d\a\wedge d\a) \}=
\nonumber \\
\frac{1}{12it} j^* K(d\bA \bA^{-1}\wedge (d\bA \bA^{-1})^2).
\end{eqnarray}
Together (\ref{do1}) and (\ref{do2}) imply the first statement of the
lemma.

2) To evaluate the form $\Otu$ on the universal vector field $v$ we
notice that
\be
da (v_{\epsilon})= - K( ad(\a) \epsilon )
\ee
for any $\epsilon \in \g$.  Taking into account (\ref{dA}) we infer
\be
\o_1(\cdot,  v_\epsilon)=\frac{1}{4it} j^* K(d\bA \bA^{-1}+\bA^{-1}
d\bA, \epsilon)
- <da, \epsilon>.
\ee
Another straightforward computation leads to
\be
\o_2(\cdot, v_{\epsilon})= \frac{1}{4it} e_{(t, u)}^* K(
A^{-1}dA-d\bar{A}\bar{A}^{-1}, A^{-1}\epsilon A- \bar{A}\epsilon
\bar{A}^{-1}).
\ee
Combining the last two equations we conclude
\be
\Otu (\cdot, v_{\epsilon})=\frac{1}{2it} e_{(t, u)}^*
K(dAA^{-1}+\bar{A}^{-1}d\bar{A}, \epsilon)- <da, \epsilon>.
\ee
Taking into account the definition (\ref{sp}) of the nondegenerate
scalar product on $\g^{C}$ one can rewrite this formula as
\be
\Otu(\cdot, v)= \frac{1}{t} e_{(t, u)}^* dA A^{-1} - da.
\ee
This observation completes the proof of Lemma 2.

{\sc Remark}

One can guess the expression (\ref{Otu}) for the 2-form $\Otu$
comparing Kirillov symplectic forms on the coadjoint orbits to the
symplectic forms on the orbits of dressing transformations computed
in \cite{GF}, \cite{AM}.

{\sc Proof of Theorem}

By assumptions of the theorem the manifold $M$ is equipped with two
maps $\m:M\rightarrow G^*_u$ and $\mu:M\rightarrow \g^*$, where $\m$
is the moment map in the sense of Lu and Weinstein and $\mu= e_{(t,
u)}^{-1}\circ \m$. Let us define a 2-form $\to$ on $M$ by the formula
\be
\to=\o - \mu^* \Otu.
\ee
In fact, the form $\to$ provides the new symplectic structure on $M$
which we are looking for.

First, observe that $\to$ is a closed 2-form on $M$:
\be
d\to=d\o-\mu^* d\Otu=0.
\ee
Moreover, $\to$ belongs to the same cohomology class as $\o$. Indeed,
$\Otu$ is a closed 2-form on the linear space $\g^*$. Hence, it is
exact. Then its pull-back $\mu^*\Otu$ is also an exact form.

Let us evaluate $\to$ on the universal vector field $v$:
\begin{eqnarray}
\to(\cdot, v)=\o(\cdot, v)-\mu^*\Otu(\cdot, v)= \nonumber \\
\frac{1}{t} \m^*(dAA^{-1})-\mu^*(\frac{1}{t} e_{(t,
u)}^*(dAA^{-1})-da)=\mu^*(da).
\end{eqnarray}
In particular, this implies that $\to$ is $G$-invariant:
\be
\L_v \to= (di_v+i_vd)\to=d\mu^*(da)=0.
\ee
So, if $\to$ defines a symplectic structure on $M$, it is
$G$-invariant and possesses an equivariant moment map
$\mu:M\rightarrow \g^*$.

The last point is to check the nondegeneracy of $\to$. Assume that at
some point $x\in M$ the form $\to$ is degenerate. This means that
there exists a nonvanishing vector $\xi$ so that
\be
\to_x(\cdot, \xi)=0.
\ee
This implies
\be
\o_x(\cdot, \xi)=\m^*  (e_{(t, u)}^{-1})^*\Otu (\cdot, \m_*\xi)\equiv
\m^* \eta.
\ee
The right hand side is a pull-back of a certain 1-form $\eta$ on
$G^*_u$ along the map $\m$.  Any such a form $\eta$ can be
represented as
\be
\eta=<dAA^{-1}, \zeta>
\ee
with some $\zeta\in \g$. Now consider a vector
\be
\tilde{\xi}=\xi-\frac{1}{t}v_{\zeta}
\ee
at the point $x\in M$. It is easy to see that  the form $\o$
annihilates this vector:
\be
\o_x(\cdot, \xi-\frac{1}{t}v_{\zeta})=\eta-t\frac{1}{t}<dAA^{-1},
\zeta>=0.
\ee
This means that the form $\o$ is also degenerate at $x$ which
contradicts to  the assumptions of the theorem. So,  $\to$ defines a
symplectic structure on $M$. This completes the proof of Theorem 1.

{\sc Remark}

It is easy to see that we can exchange the roles of Hamiltonian and
Poisson actions in Theorem 1 . Moreover, we can directly compare
Poisson actions with different values of parameters $t$ and $u$.

\subsection*{Corollaries for Poisson actions}

Here we give new short proves of several results on the actions of
Poisson-Lie groups on symplectic manifolds.

Recently Flashka and Ratiu  \cite{FR} proved the following convexity
theorem for the moment map in the sense of Lu and Weinstein (see also
\cite{K}, \cite{LR}).

\begin{cor}
Let $M$ be a compact symplectic manifold which carries a Poisson
action $\A$ of the compact  group $G$ equipped with the Poisson
structure $\P_{(t, u)}$. Assume that there exists an equivariant
moment map $\m:M\rightarrow G^*_u$. Define the map $\mu=e_{(t,
u)}^{-1}\circ \m$. Then the intersection of $\mu(M)$ with the
positive Weyl chamber $W_+$
\be
\mu_+(M)=\mu(M)\cap W_+
\ee
is a convex polytop.
\end{cor}

{\sc Proof}

As we know, the map $\mu$ provides a Hamiltonian equivariant moment
map for some symplectic structure on $M$. Convexity property for the
map $\m$ as stated above coincides with the standard convexity for
the Hamiltonian moment map $\mu$ \cite{A}, \cite{GS}.

The technique of Hamiltonian reduction has been generalized to
Poisson actions by Lu \cite{Lu}.
Here we need some new notations and definitions to formulate a
statement.

\begin{defn}
The value $\gamma\in G^*_u$ is called a regular value of the moment
map $\m:M\rightarrow G^*_u$ if  some  quotient of $G$ over a discrete
 (possibly trivial) subgroup $F$ of the center of $G$ acts freely on
$\m^{-1}(\gamma)$.
\end {defn}
It is convenient to introduce a special notation for the canonical
projection
\be
\pi:M\rightarrow M/G
\ee
to the quotient space $M/G$ and for embedding of   $\m^{-1}(\gamma)$
into $M$:
\be
i_{\gamma}:\m^{-1}(\gamma) \rightarrow  M.
\ee

\begin{cor}
Let $M$ be a  symplectic manifold which carries a Poisson action $\A$
of the compact  group $G$ equipped with the Poisson structure
$\P_{(t, u)}$. Assume that there exists an equivariant moment map
$\m:M\rightarrow G^*_u$. Choose some $\gamma\in G^*_u$ being a
regular value of the moment map.  Then
$M_{\gamma}=\pi(\m^{-1}(\gamma))$ is a symplectic manifold with
symplectic structure $\o_{\gamma}$ defined via
\be
\pi^*\o_{\gamma}=i^*_{\gamma} \o.
\ee
\end{cor}

{\sc Proof}

Let us switch to the symplectic structure $\to$ on $M$ and let
$c=e_{(t, u)}^{-1}(\gamma)$. The map $e_{(t, u)}$ being equivariant,
the space $M_{\gamma}$ coincides with the reduced space
obtained by the Hamiltonian reduction over the value $c$ of the
moment map $\mu$. In fact, symplectic structures of the Hamiltonian
and Poisson reduced spaces  also coincide as
\be
i_{\gamma}^*(\o-\to)=i^*_{\gamma}\mu^*\Otu=0.
\ee
The latter is true because the embedding $i_{\gamma}$ chooses the
point in $c\in \g^*$ and the pull-back of the 2-form $\Otu$ to this
point vanishes for dimensional reasons.

By now we compared $(M, \o)$ and $(M, \to)$ as symplectic $G$-spaces.
It is clear that they do not coincide in this category as the
$G$-action preserves $\to$ and changes $\o$. However, it possible
that $(M, \o)$ and $(M, \to)$ are isomorphic as symplectic spaces
(now we disregard the $G$-action).  This is indeed the case, the
isomorphism between $(M, \o)$ and $(M, \to)$ is called
Ginzburg-Weinstein isomorphism \cite{GW}.

\begin{cor}
For arbitrary values of parameters $t$ and $u$ $(M, \o)$ and $(M,
\to)$ are isomorphic as symplectic spaces. In particular, orbits of
dressing transformations are symplectomorphic to the corresponding
coadjoint orbits.
\end{cor}

{\sc Proof}

Choose some primitive $\a_{(t, u)}$ of the 2-form $\Otu$:
\be
\Otu=d \a_{(t, u)}.
\ee
We would like to vary parameters $t$ and $u$ of the Poisson bracket
of $G$. For simplicity we change only  $t$. When $t$ varies, the form
the symplectic form $\o=\to+\mu^*\Otu$ changes as:
\be
\frac{\partial}{\partial t} \o= \mu^* \ \frac{\partial \Otu}{\partial
t}= \mu^* d \ \frac{\partial \a_{(t, u)}}{\partial t}.
\ee
Denote
\be
\beta_{(t, u)}=\frac{\partial \a_{(t, u)}}{\partial t}
\ee
and construct a vector field $V_{(t, u)}$
\be
V_{(t, u)}=\P_M(\cdot,  \mu^*\beta_{(t, u)}).
\ee
The vector field $V_{(t, u)}$ is a certain linear combination of the
vector fields $v_{\epsilon}$ with coefficients which depend only on
the value of the moment map $\m(x)$:
\be \label{v}
V_{(t, u)}=<{\cal E}(\m(x)), v>\ , \ \beta_{(t, u)}=<{\cal E}(A),
dAA^{-1}>.
\ee
The Lie derivative of the symplectic structure $\o$ with respect to
$V_{(t, u)}$ coincides
with the $t$-derivative:
\be
\L_{V_{(t, u)}}\o =di_{V_{(t, u)}}\o=d \mu^*\beta_{(t, u)}=
\frac{\partial \o}{\partial t}.
\ee
Integrating the ($t$-dependent) field $V_{(t, u)}$ we construct a
family of Ginzburg-Weinstein isomorphisms which identify $(M, \o)$
and $(M, \to)$ for different values of $t$. One can
construct symplectomorphisms between these spaces with different
values of $u$ in a similar fashion.

{\sc Remark}

Formula (\ref{v}) for the  vector field $V_{(t, u)}$ makes it
possible to extend the Ginzburg-Weinstein isomorphism to Poisson
manifolds  which carry a Poisson $G$-action and possess an
equivariant moment map $\m:M\rightarrow G^*_u$ in the following
sense:
\be
v=\frac{1}{t} \P_M(\cdot,  \m^*( dAA^{-1})).
\ee
This condition  implies that symplectic leaves are preserved by the
$G$-action.
Integrating the vector field (\ref{v}) one can obtain a
diffeomorphism $D_{(t, u)}$ of $M$ which
preserves symplectic leaves  and replaces the Poisson structure
$\P_M$ by the
$G$-invariant Poisson structure $\tilde{\P}_M$. Restricted to each
symplectic leaf
$D_{(t, u)}$ coincides with the Ginzburg-Weinstein symplectomorphism
described
above.  This implies that the new Poisson $G$-space $(M,
\tilde{\P}_M)$ possesses an equivariant moment map $\mu=e_{(t,
u)}^{-1}\circ \m$  which arises from the equivariant
moment maps on each symplectic leaf.

Let us apply this construction to the Poisson space $G^*_u$ equipped
with the Poisson structure $\P^*_{(t, u)}$ from the standard list.
The dressing action of $G$
preserves symplectic leaves, the moment map is equal to identity
$\m=id:G^*_u\rightarrow G^*_u$.  The Ginzburg-Weinstein
diffeomorphism $D_{(t, u)}$ endows $G^*_u$ with a new $G$-invariant
Poisson structure $\tilde{\P}^*_{(t, u)}$ and a new moment map
$\mu=e_{(t, u)}^{-1}:G^*_u\rightarrow \g^*$.  Both maps $D_{(t, u)}$
and $\mu$ are invertible Poisson maps. Thus, an invertible Poisson
map $e_{(t, u)}^{-1}\circ D_{(t, u)}$ establishes a Poisson
isomorphism of $(G^*_u, \P^*_{(t, u)})$  and $\g^*$ equipped with the
standard Kirillov-Kostant-Sourieu bracket.   In fact, we have
recovered the original version of  the Ginzburg-Weinstein isomorphism
\cite{GW}.

\subsection*{Appendix. The case of t=0}

Here we collect some details on the special family of Poisson
structures on compact Lie groups  which may be obtained from the
general case (\ref{fam}) in the limit
\be \label{spec}
(t\rightarrow 0, u=\frac{w}{t})\ , w=const.
\ee

All results obtained in the main text generalize to the special
family (\ref{spec}). In fact, in this limit calculations become much
easier. For this reason, we provide only the basic definitions and
formulas related to the proof of Lemma 2. The proves of Theorem 1 and
of all Corollaries do
not change.

For the special family of Poisson structures (\ref{spec}) the dual
Lie algebra is a subset in the semi-direct product of the Cartan
subalgebra $\h$ and the the dual Lie algebra $\g^*_0$ considered as
an Abelian Lie algebra:
\be \label{g0w}
\g^*_{(0, w)}=\{ (ih+n, -w(h)), h\in \h, n\in \n^C\}.
\ee
The $\h$ component acts on the $\g^*_0$ component by the natural
coadjoint action.

The corresponding Lie group is a subgroup in the semi-direct product
of the Cartan subgroup $H$ and $\g^*_0$  (viewed  as an Abelian group
with addition playing the role of the group operation):
\be
G^*_{(0, w)}=\{ (ih+n, exp\{-w(h)\} ), h\in \h, n\in \n^C\}.
\ee

The equivariant map $e_w:\g^*_0\rightarrow G^*_{(0, w)}$ is defined
as
\be
e_w(ih+n)=(ih+n, exp\{-w(h)\} ).
\ee
The inverse map $e_w^{-1}$ is a forgetting map which drops the second
component of the pair.

It is instructive to compare Maurer-Cartan forms for the Abelian
group  $\g^*_0$:
\be
a=ih+n\ , \ da= idh+dn
\ee
and for the group $G^*_{(0, w)}$:
\begin{eqnarray}
& A=(ih+n, exp\{ -w(h)\} )& \nonumber \\
& dAA^{-1}=(idh+dn-[w(dh),n], -w(dh)). &
\end{eqnarray}
Let us mention that the second component in the pair describing
$dAA^{-1}$ is disregarded in the pairing with elements of $\g$.

The definition of the  moment map in the sense of Lu and Weinstein
modifies as follows:
\begin{defn}
Let $G$ be a compact Poisson-Lie group equipped with a Poisson
structure $\P_{(0, w)}$. Let $\A: G\times M \rightarrow M$ be a
Poisson action of $G$ on the symplectic manifold $M$
The map $\m: M\rightarrow G^*_{(0, w)}$  is called a moment map in
the sense of Lu and Weinstein if
\be
\o(\cdot, v)= \m^* (dAA^{-1}),
\ee
where  $dAA^{-1}$ is a right-invariant Maurer-Cartan form on
$G^*_{(0, w)}$.
\end{defn}

Lemma 2 in this situation is reformulated as:

\begin{lemma}
There exists such a 2-form $\Omega_w$ on $\g^*$, so that the
following two properties are fulfilled:

\vskip 0.2cm
1) The form $\Omega_w$ is closed $d \Omega_w =0$.

\vskip 0.2cm
2) $\Omega_w(\cdot, v)= e_w^* dA A^{-1} - da$.

\vskip 0.2cm
Here $v$ is   the universal vector field corresponding to the
coadjoint action of $G$ on $\g^*$,
$a\in g^*$ and  $A=e_w (a) \in G^*_{(0, w)}$.
\end{lemma}

{\sc Proof}

The 2-form $\Omega_w$ which fulfils these two properties looks as
\be
\Omega_w=\frac{1}{2} w(dh\wedge dh),
\ee
where $h$ is a Cartan projection of $(ih+n)\in \g^*_0$.

Obviously,  $\Omega_w$ is closed. Evaluating  it on the universal
vector field $v$ one finds:
\begin{eqnarray} \label{2l2}
\Omega_w(\cdot, v_{\epsilon}) &= & \frac{1}{2i} K(w(dh), [\epsilon,
n+\bar{n}])=\nonumber \\
-<[w(dh), n], \epsilon> & = & <e_w^*dAA^{-1}-da, \epsilon>
\end{eqnarray}
This completes the proof of Lemma 3.

\subsection*{Acknowledgements}

This work has been done during the Workshop on Geometry and Physics
held in Aahrus in Summer 1995. I have been inspired by the work of
Flaschka and Ratiu \cite{FR} on the convexity theorem for the
Poisson-Lie moment map. I am indebted to P.Cartier, Y.Karshon and
V.Rubtsov  for discussions and comments.

\end{document}